\documentclass[11pt]{article}
\usepackage{cite}
\usepackage{mathrsfs}
\usepackage{amsfonts}
\usepackage{amsmath}
\usepackage{amsfonts,amssymb}
\usepackage{dsfont}
\usepackage{curves}
\usepackage{mathrsfs}
\usepackage{pifont}
\usepackage{amssymb}
\allowdisplaybreaks

\numberwithin{equation}{section}

\date{}

\def\BigRoman{\uppercase\expandafter{\romannumeral\number\count 255 }}
\def\Romannumeral{\afterassignment\BigRoman\count255=}

\begin{document}
\title{A result on spanning trees with bounded total excess
}
\author{\small  Sizhong Zhou\footnote{Corresponding
author. E-mail address: zsz\_cumt@163.com (S. Zhou)}\\
\small  School of Science, Jiangsu University of Science and Technology,\\
\small  Zhenjiang, Jiangsu 212100, China\\
}

\maketitle
\begin{abstract}
\noindent Let $G$ be a connected graph and $T$ a spanning tree of $G$. Let $\rho(G)$ denote the adjacency spectral radius of $G$. The $k$-excess
of a vertex $v$ in $T$ is defined as $\max\{0,d_T(v)-k\}$. The total $k$-excess $\mbox{te}(T,k)$ is defined by $\mbox{te}(T,k)=\sum\limits_{v\in V(T)}{\max\{0,d_T(v)-k\}}$.
A tree $T$ is said to be a $k$-tree if $d_T(v)\leq k$ for any $v\in V(T)$, that is to say, the maximum degree of a
$k$-tree is at most $k$. In fact, $T$ is a spanning $k$-tree if and only if $\mbox{te}(T,k)=0$. This paper studies a generalization of spanning
$k$-trees using a concept called total $k$-excess and proposes a lower bound for $\rho(G)$ in a connected graph $G$ to ensure that $G$ contains
a spanning tree $T$ with $\mbox{te}(T,k)\leq b$, where $k$ and $b$ are two nonnegative integers with $k\geq\max\{5,b+3\}$ and $(b,k)\neq(2,5)$.
\\
\begin{flushleft}

(2020) Mathematics Subject Classification: 05C50, 05C05
\end{flushleft}
\end{abstract}

\section{Introduction}

Throughout this paper, we only consider finite and undirected graphs without loops or multiple edges. Let $G$ denote a graph. Use $V(G)$ and
$E(G)$ to denote the vertex set and the edge set of $G$, respectively. For a vertex $v$ of $G$, we write the degree of $v$ in $G$ by $d_G(v)$.
A vertex $v$ of $G$ with $d_G(v)=0$ is called an isolated vertex of $G$. For a given subset $S$ of $V(G)$, we let $|S|$ and $G[S]$ denote the
number of elements in $S$ and the subgraphs of $G$ induced by $S$, respectively. We write $G-S$ for $G[V(G)\setminus S]$. The number of
components in $G-S$ is denoted by $c(G-S)$. Let $K_n$ denote the complete graph of order $n$. For two given graphs $G_1$ and $G_2$, the union
of $G_1$ and $G_2$ is denoted by $G_1\cup G_2$. The join $G_1\vee G_2$ is obtained from $G_1\cup G_2$ by joining every vertex of $G_1$ with
every vertex of $G_2$ by an edge.

For an integer $k\geq2$ and a connected graph $G$, a spanning $k$-tree $T$ of $G$ is a spanning tree whose maximum degree is at most $k$. In
particular, a spanning 2-tree is also called a Hamilton path. Obviously, a spanning $k$-tree of $G$ is a natural generalization of a Hamilton
path. Enomoto, Ohnishi and Ota \cite{EOO} introduced the notion of total excess. For a spanning tree $T$ of $G$, the $k$-excess of a vertex $v$
in $T$ is defined as $\max\{0,d_T(v)-k\}$. The total $k$-excess $\mbox{te}(T,k)$ is defined by $\mbox{te}(T,k)=\sum\limits_{v\in V(T)}{\max\{0,d_T(v)-k\}}$.
In fact, $T$ is a spanning $k$-tree if and only if $\mbox{te}(T,k)=0$.

Win \cite{Win} verified the following result, which provides a sufficient condition for a connected graph to have a spanning $k$-tree. Ellingham
and Zha \cite{EZ} gave a short proof.

\medskip

\noindent{\textbf{Theorem 1.1}} (Win \cite{Win}, Ellingham and Zha \cite{EZ}). Let $k$ be an integers with $k\geq2$, and let $G$ be a connected graph.
If
$$
c(G-S)\leq(k-2)|S|+2
$$
for any vertex subset $S$ of $G$, then $G$ has a spanning tree.

\medskip

Enomoto, Ohnishi and Ota \cite{EOO} gave a generalization of Theorem 1.1 based on total $k$-excess $\mbox{te}(T,k)$.

\medskip

\noindent{\textbf{Theorem 1.2}} (Enomoto, Ohnishi and Ota \cite{EOO}). Let $b$ and $k$ be two nonnegative integers with $k\geq2$, and let $G$ be a
connected graph. If for any vertex subset $S$ of $G$
$$
c(G-S)\leq(k-2)|S|+b+2,
$$
then $G$ has a spanning tree $T$ with $\mbox{te}(T,k)\leq b$.

\medskip

Ota and Sugiyama \cite{OS} showed a sufficient condition using the condition on forbidden subgraphs for a graph to contain a spanning $k$-tree.
Rose \cite{Rose}
obtained some characterizations for graphs to possess spanning $k$-trees. Furuya et al \cite{FMMMTY} claimed a degree sum condition for the
existence of a spanning $k$-tree in a connected $K_{1,k+1}$-free graph. Ohnishi, Ota and Ozeki \cite{OOO} provided some sufficient conditions
for a graph to contain a spanning tree with bounded total excess. Ozeki \cite{Ozeki} presented a toughness condition for a graph to have a
spanning tree with bounded total excess. Maezawa, Tsugaki and Yashima \cite{MTY} gave a characterization for a graph to have a spanning tree
with bounded total excess.

Given a graph $G$ with $V(G)=\{v_1,v_2,\ldots,v_n\}$, the adjacency matrix of $G$ is the $n\times n$ matrix whose rows and columns are indexed
by vertices and is defined by $A(G)=(a_{ij})_{n\times n}$, where
\[
a_{ij}=\left\{
\begin{array}{ll}
1,&if \ v_iv_j\in E(G);\\
0,&otherwise.\\
\end{array}
\right.
\]
The eigenvalues of $A(G)$ are called the adjacency eigenvalues of $G$, which are denoted by $\lambda_1(G),\lambda_2(G),\ldots,\lambda_n(G)$ and
are arranged as $\lambda_1(G)\geq\lambda_2(G)\geq\cdots\geq\lambda_n(G)$. Notice that the adjacency spectral radius of $G$ equals $\lambda_1(G)$,
written as $\rho(G)$. More information on the spectral radius can be found in \cite{O,FN,LSY,Ws,Ws1,AN,Zs1,ZZS,Zt,ZZL2,Zs,ZZL1,ZZ}.

Gu and Liu \cite{GL} provided a Laplacian eigenvalue condition for a connected graph to contain a spanning $k$-tree. Wu \cite{Wc} established
a lower bound on the spectral radius of a connected graph $G$ to ensure that $G$ has a spanning tree with leaf degree at most $k$. Zhou, Sun and
Liu \cite{ZSL} obtained a distance spectral radius and a distance signless Laplacian spectral radius condition for the existence of spanning
trees with leaf degree at most $k$ in graphs, respectively. Fan et al \cite{FGHL} proposed an adjacency spectral radius and a signless Laplacian
spectral radius condition for a connected graph to contain a spanning $k$-tree, respectively. Zhou and Wu \cite{ZW} proved an upper bound for
the distance spectral radius of a connected graph $G$ to ensure the existence of a spanning $k$-tree in $G$. Zhou, Zhang and Liu \cite{ZZL}
showed a distance signless Laplacian spectral radius condition for a connected graph to contain a spanning $k$-tree.

Ning and Ge \cite{NG} proposed a sufficient condition in terms of the adjacency spectral radius for the existence of a spanning 2-tree (or
Hamiltonian path) in a connected graph.

\medskip

\noindent{\textbf{Theorem 1.3}} (Ning and Ge \cite{NG}). Let $G$ be a connected graph of order $n\geq4$. If
$$
\rho(G)>n-3,
$$
then $G$ has a spanning 2-tree (or Hamiltonian path), unless $G\in\{K_1\vee(K_{n-3}\cup2K_1,K_2\vee4K_1,K_1\vee(K_{1,3}\cup K_1)\}$.

\medskip

As a generalization of Theorem 1.3, Fan et al \cite{FGHL} provided a sufficient condition based on the adjacency spectral radius for a connected
graph to possess a spanning $k$-tree.

\medskip

\noindent{\textbf{Theorem 1.4}} (Fan et al \cite{FGHL}). Let $G$ be a connected graph of order $n\geq2k+16$, where $k$ is an integer with $k\geq3$.
If
$$
\rho(G)\geq\rho(K_1\vee(K_{n-k-1}\cup kK_1)),
$$
then $G$ contains a spanning $k$-tree, unless $G=K_1\vee(K_{n-k-1}\cup kK_1)$.

\medskip

Motivated by \cite{EOO,FGHL} directly, we investigate the existence of spanning trees with bounded total excess in connected graphs and obtain an
adjacency spectral condition for connected graphs to contains spanning trees with bounded total excess. Our main result is a generalization of
Theorem 1.4 based on total $k$-excess $\mbox{te}(T,k)$.

\medskip

\noindent{\textbf{Theorem 1.5.}} Let $k$ and $b$ be two nonnegative integers with $k\geq\max\{5,b+3\}$ and $(b,k)\neq(2,5)$, and let $G$ be a
connected graph of order $n\geq k+b+2$. If
$$
\rho(G)\geq\rho(K_1\vee(K_{n-k-b-1}\cup(k+b)K_1)),
$$
then $G$ contains a spanning tree $T$ with $\mbox{te}(T,k)\leq b$, unless $G=K_1\vee(K_{n-k-b-1}\cup(k+b)K_1)$.

\medskip

We immediately obtain the following corollary from Theorem 1.5 by setting $b=0$.

\medskip

\noindent{\textbf{Corollary 1.6.}} Let $k$ be an integers with $k\geq5$, and let $G$ be a connected graph of order $n\geq k+2$. If
$$
\rho(G)\geq\rho(K_1\vee(K_{n-k-1}\cup kK_1)),
$$
then $G$ contains a spanning $k$-tree, unless $G=K_1\vee(K_{n-k-1}\cup kK_1)$.

\medskip

When $k\geq5$, the order of Corollary 1.6 is better than that of Theorem 1.4, and our improvement on the order in Corollary 1.6 is very large.
Moreover, Corollary 1.6 is immediately obtained from Theorem 1.5. From the above discussion, Theorem 1.5 is a generalization and an improvement
of Theorem 1.4. The paper extends the existing body of knowledge on spanning $k$-trees by focusing on the spectral radius, which can enrich the
spectral graph theory and its applications in studying the structural properties of graphs. Consequently, the result of Theorem 1.5 is interesting
and important.

\section{Preliminary lemmas}

In this section, we list some lemmas which will be used in the proof of our main result.

\medskip

\noindent{\textbf{Lemma 2.1}} (Brouwer and Haemers \cite{BH}). Let $G$ be a connected graph, and let $H$ be a subgraph of $G$. Then
$$
\rho(G)\geq\rho(H),
$$
with equality holding if and only if $H=G$.

\medskip

\noindent{\textbf{Lemma 2.2}} (Fan et al \cite{FGHL}). Let $n=\sum\limits_{i=1}^{t}n_i+s$. If $n_1\geq n_2\geq\cdots\geq n_t\geq1$ and $n_1<n-s-t+1$,
then
$$
\rho(K_s\vee(K_{n_1}\cup K_{n_2}\cup\cdots\cup K_{n_t}))<\rho(K_s\vee(K_{n-s-t+1}\cup(t-1)K_1)).
$$

\medskip

Let $M$ be a real $n\times n$ matrix described in the following block form
\begin{align*}
M=\left(
  \begin{array}{cccc}
    M_{11} & M_{12} & \cdots & M_{1m}\\
    M_{21} & M_{22} & \cdots & M_{2m}\\
    \vdots & \vdots & \ddots & \vdots\\
    M_{m1} & M_{m2} & \cdots & M_{mm}\\
  \end{array}
\right),
\end{align*}
whose rows and columns are partitioned into subsets $X_1,X_2,\ldots,X_m$ of $\{1,2,\ldots,n\}$. Let $M_{ij}$ denote the block of $M$ by deleting
the rows in $\{1,2,\ldots,n\}-X_i$ and the columns in $\{1,2,\ldots,n\}-X_j$. A quotient matrix $B$ is the $m\times m$ matrix whose entries are
the average row sums of the blocks $M_{ij}$ of $M$. The partition is called equitable if every block $M_{ij}$ of $M$ has a constant row sum.

\medskip

\noindent{\textbf{Lemma 2.3}} (You, Yang, So and Xi \cite{YYSX}). Let $M$ be a real symmetric matrix with an equitable partition and let $B$ be
the corresponding quotient matrix. Then the eigenvalues of $B$ are also eigenvalues of $M$. Furthermore, if $M$ is a nonnegative matrix, then
the spectral radius of $B$ is equal to the spectral radius of $M$.

\section{The proof of Theorem 1.5}

\noindent{\it Proof of Theorem 1.5.} Suppose to the contrary that $G$ contains no spanning tree $T$ with $\mbox{te}(T,k)\leq b$. According to
Theorem 1.2, there exists some nonempty subset $S$ of $V(G)$ such that $c(G-S)\geq(k-2)|S|+b+3$. For convenience, we let $|S|=s$ and $c(G-S)=c$.
Then $c\geq(k-2)s+b+3$, $n\geq(k-1)s+b+3$ and $G$ is a spanning subgraph of $G_1=K_s\vee(K_{n_1}\cup K_{n_2}\cup\cdots\cup K_{n_{(k-2)s+b+3}})$,
where $n_1,n_2,\ldots,n_{(k-2)s+b+3}$ are positive integers with $\sum\limits_{i=1}^{(k-2)s+b+3}{n_i}=n-s$. Without loss of generality, we may
assume $n_1\geq n_2\geq\cdots\geq n_{(k-2)s+b+3}\geq1$. In terms of Lemma 2.1, we infer
\begin{align}\label{eq:3.1}
\rho(G)\leq\rho(G_1),
\end{align}
where the equality holds if and only if $G=G_1$. Let $G_2=K_s\vee(K_{n-(k-1)s-b-2}\cup((k-2)s+b+2)K_1)$. Using Lemma 2.2, we conclude
\begin{align}\label{eq:3.2}
\rho(G_1)\leq\rho(G_2),
\end{align}
where the equality occurs if and only if $G_1=G_2$. By virtue of the partition $V(G_2)=V(K_s)\cup V(K_{n-(k-1)s-b-2})\cup V(((k-2)s+b+2)K_1)$,
the quotient matrix of $A(G_2)$ is
\begin{align*}
B_1=\left(
  \begin{array}{ccc}
    s-1 & n-(k-1)s-b-2 & (k-2)s+b+2\\
    s & n-(k-1)s-b-3 & 0\\
    s & 0 & 0\\
  \end{array}
\right).
\end{align*}
By a simple computation, the characteristic polynomial of $B_1$ equals
\begin{align*}
\varphi_{B_1}(x)=&x^{3}+(-n+(k-2)s+b+4)x^{2}-(n+(k-2)s^{2}-(k-b-4)s-b-3)x\\
&+(k-2)s^{2}n+(b+2)sn-(k-2)(k-1)s^{3}-(2bk+5k-3b-8)s^{2}\\
&-(b+2)(b+3)s.
\end{align*}
Notice that the partition $V(G_2)=V(K_s)\cup V(K_{n-(k-1)s-b-2})\cup V(((k-2)s+b+2)K_1)$ is equitable. Using Lemma 2.3, the largest root, say $\rho_1$,
of $\varphi_{B_1}(x)=0$ satisfies $\rho_1=\rho(G_2)$. Note that $K_s\vee((k-2)s+b+3)K_1$ is a subgraph of $G_2$. In view of the partition
$V(K_s\vee((k-2)s+b+3)K_1)=V(K_s)\cup V(((k-2)s+b+3)K_1)$, the quotient matrix of $A(K_s\vee((k-2)s+b+3)K_1)$ is
\begin{align*}
B_2=\left(
  \begin{array}{ccc}
    s-1 & (k-2)s+b+3\\
    s & 0\\
  \end{array}
\right).
\end{align*}
Then the characteristic polynomial of $B_2$ is
$$
\varphi_{B_2}(x)=x^{2}-(s-1)x-s((k-2)s+b+3).
$$
Since the partition $V(K_s\vee((k-2)s+b+3)K_1)=V(K_s)\cup V(((k-2)s+b+3)K_1)$ is equitable, it follows from Lemma 2.3 that $\rho(K_s\vee((k-2)s+b+3)K_1)$
is the largest root of $\varphi_{B_2}(x)=0$. By a direct calculation, we obtain
\begin{align}\label{eq:3.3}
\rho(K_s\vee((k-2)s+b+3)K_1)=\frac{s-1+\sqrt{(4k-7)s^{2}+(4b+10)s+1}}{2}.
\end{align}
According to \eqref{eq:3.3} and Lemma 2.1, we deduce
\begin{align}\label{eq:3.4}
\rho_1=&\rho(G_2)\nonumber\\
\geq&\rho(K_s\vee((k-2)s+b+3)K_1)\nonumber\\
=&\frac{s-1+\sqrt{(4k-7)s^{2}+(4b+10)s+1}}{2}.
\end{align}

Let $G_*=K_1\vee(K_{n-k-b-1}\cup(k+b)K_1)$. Then the quotient matrix of $A(G_*)$ in view of the partition $V(G_*)=V(K_1)\cup V(K_{n-k-b-1})\cup V((k+b)K_1)$
can be written as
\begin{align*}
B_*=\left(
  \begin{array}{ccc}
    0 & n-k-b-1 & k+b\\
    1 & n-k-b-2 & 0\\
    1 & 0 & 0\\
  \end{array}
\right).
\end{align*}
By a simple computation, the characteristic polynomial of $B_*$ equals
\begin{align*}
\varphi_{B_*}(x)=x^{3}+(-n+k+b+2)x^{2}-(n-1)x+(k+b)n-k^{2}-2bk-2k-b^{2}-2b.
\end{align*}
Notice that the partition $V(G_*)=V(K_1)\cup V(K_{n-k-b-1})\cup V((k+b)K_1)$ is equitable. From Lemma 2.3, $\rho(G_*)$ is the largest root of
$\varphi_{B_*}(x)=0$.

For $s=1$, we possess $G_2=K_1\vee(K_{n-k-b-1}\cup(k+b)K_1)$, and so $\rho(G_2)=\rho(K_1\vee(K_{n-k-b-1}\cup(k+b)K_1))$. Combining this with \eqref{eq:3.1}
and \eqref{eq:3.2}, we admit
$$
\rho(G)\leq\rho(K_1\vee(K_{n-k-b-1}\cup(k+b)K_1)),
$$
with equality holding if and only if $G=K_1\vee(K_{n-k-b-1}\cup(k+b)K_1)$. Together with the condition $\rho(G)\geq\rho(K_1\vee(K_{n-k-b-1}\cup(k+b)K_1))$ of
the theorem, we claim $G=K_1\vee(K_{n-k-b-1}\cup(k+b)K_1)$, which is a contradiction to $G\neq K_1\vee(K_{n-k-b-1}\cup(k+b)K_1)$. In what follows, we consider
$s\geq2$.

By plugging the value $\rho_1$ into $x$ of $\varphi_{B_*}(x)-\varphi_{B_1}(x)$, we obtain
\begin{align}\label{eq:3.5}
\varphi_{B_*}(\rho_1)=\varphi_{B_*}(\rho_1)-\varphi_{B_1}(\rho_1)=(s-1)f_1(\rho_1),
\end{align}
where $f_1(\rho_1)=-(k-2)\rho_1^{2}+((k-2)s+b+2)\rho_1-(k-2)sn-(k+b)n+(k-2)(k-1)s^{2}+(k^{2}+2bk+2k-3b-6)s+k^{2}+2bk+b^{2}+2k+2b$. Since
$n\geq(k-1)s+b+3$, we get
\begin{align}\label{eq:3.6}
f_1(\rho_1)=&-(k-2)\rho_1^{2}+((k-2)s+b+2)\rho_1-(k-2)sn-(k+b)n\nonumber\\
&+(k-2)(k-1)s^{2}+(k^{2}+2bk+2k-3b-6)s+k^{2}+2bk+b^{2}+2k+2b\nonumber\\
\leq&-(k-2)\rho_1^{2}+((k-2)s+b+2)\rho_1-(k-2)s((k-1)s+b+3)\nonumber\\
&-(k+b)((k-1)s+b+3)+(k-2)(k-1)s^{2}+(k^{2}+2bk+2k-3b-6)s\nonumber\\
&+k^{2}+2bk+b^{2}+2k+2b\nonumber\\
=&-(k-2)\rho_1^{2}+((k-2)s+b+2)\rho_1+k^{2}+bk-k-b\nonumber\\
\triangleq&f_2(\rho_1).
\end{align}
The symmetry axis of $f_2(\rho_1)$ is $\rho_1=\frac{(k-2)s+b+2}{2(k-2)}$, which implies that $f_2(\rho_1)$ is decreasing for $\rho_1\geq\frac{(k-2)s+b+2}{2(k-2)}$.

By $s\geq2$ and $k\geq\max\{5,b+3\}$, we have
\begin{align*}
&(k-2)^{2}((4k-7)s^{2}+(4b+10)s+1)-((k-2)s+b+2-(k-2)(s-1))^{2}\\
=&(k-2)^{2}((4k-7)s^{2}+(4b+10)s+1)-(b+k)^{2}\\
\geq&(k-2)^{2}(4(4k-7)+2(4b+10)+1)-(b+k)^{2}\\
=&(k^{2}-4k+4)(16k+8b-7)-b^{2}-2bk-k^{2}\\
\geq&(k+4)(16k+8b-7)-b^{2}-2bk-k^{2}\\
=&15k^{2}+6bk+57k-b^{2}+32b-28\\
\geq&15(b+3)^{2}+6b(b+3)+57(b+3)-b^{2}+32b-28\\
=&20b^{2}+197b+278\\
>&0,
\end{align*}
which implies
$$
\frac{(k-2)s+b+2}{2(k-2)}<\frac{s-1+\sqrt{(4k-7)s^{2}+(4b+10)s+1}}{2}.
$$
Together with \eqref{eq:3.4}, we get
$$
\frac{(k-2)s+b+2}{2(k-2)}<\frac{s-1+\sqrt{(4k-7)s^{2}+(4b+10)s+1}}{2}\leq\rho_1.
$$
Recall that $f_2(\rho_1)$ is decreasing for $\rho_1\geq\frac{(k-2)s+b+2}{2(k-2)}$. Combining these with \eqref{eq:3.6}, we obtain
\begin{align}\label{eq:3.7}
f_1(\rho_1)\leq&f_2(\rho_1)\nonumber\\
\leq&f_2\left(\frac{s-1+\sqrt{(4k-7)s^{2}+(4b+10)s+1}}{2}\right)\nonumber\\
=&-(k-2)\left(\frac{s-1+\sqrt{(4k-7)s^{2}+(4b+10)s+1}}{2}\right)^{2}\nonumber\\
&+((k-2)s+b+2)\left(\frac{s-1+\sqrt{(4k-7)s^{2}+(4b+10)s+1}}{2}\right)\nonumber\\
&+k^{2}+bk-k-b\nonumber\\
=&\frac{1}{2}(-(2k^{2}-8k+8)s^{2}-(2bk+5k-5b-12)s+2k^{2}+2bk-3k-3b\nonumber\\
&+(k+b)\sqrt{(4k-7)s^{2}+(4b+10)s+1}).
\end{align}

Notice that
\begin{align*}
&k^{2}s^{2}-((4k-7)s^{2}+(4b+10)s+1)\\
&=(k^{2}-4k+7)s^{2}-(4b+10)s-1\\
&\geq(2k^{2}-8k+14)s-(4b+10)s-1\\
&=(2k^{2}-8k-4b+4)s-1\\
&\geq(2k^{2}-8k-4(k-3)+4)s-1\\
&=(k-4)(2k-4)s-1\\
&>0
\end{align*}
by $s\geq2$ and $k\geq\max\{5,b+3\}$. Combining this with \eqref{eq:3.7}, $s\geq2$, $k\geq\max\{5,b+3\}$ and $(b,k)\neq(2,5)$, we obtain
\begin{align}\label{eq:3.8}
f_1(\rho_1)<&\frac{1}{2}(-(2k^{2}-8k+8)s^{2}-(2bk+5k-5b-12)s\nonumber\\
&+2k^{2}+2bk-3k-3b+(k+b)ks)\nonumber\\
=&\frac{1}{2}(-(2k^{2}-8k+8)s^{2}+(k^{2}-bk-5k+5b+12)s\nonumber\\
&+2k^{2}+2bk-3k-3b)\nonumber\\
\leq&\frac{1}{2}(-(4k^{2}-16k+16)s+(k^{2}-bk-5k+5b+12)s\nonumber\\
&+2k^{2}+2bk-3k-3b)\nonumber\\
=&\frac{1}{2}((-3k^{2}-bk+11k+5b-4)s+2k^{2}+2bk-3k-3b)\nonumber\\
\leq&\frac{1}{2}(-6k^{2}-2bk+22k+10b-8+2k^{2}+2bk-3k-3b)\nonumber\\
=&\frac{1}{2}(-4k^{2}+19k+7b-8)\nonumber\\
<&0.
\end{align}

It follows from \eqref{eq:3.5}, \eqref{eq:3.8} and $s\geq2$ that
$$
\varphi_{B_*}(\rho_1)=(s-1)f_1(\rho_1)<0,
$$
which implies
\begin{align}\label{eq:3.9}
\rho(G_2)=\rho_1<\rho(G_*)=\rho(K_1\vee(K_{n-k-b-1}\cup(k+b)K_1)).
\end{align}
According to \eqref{eq:3.1}, \eqref{eq:3.2} and \eqref{eq:3.9}, we deduce
$$
\rho(G)\leq\rho(G_1)\leq\rho(G_2)<\rho(K_1\vee(K_{n-k-b-1}\cup(k+b)K_1)),
$$
which contradicts $\rho(G)\geq\rho(K_1\vee(K_{n-k-b-1}\cup(k+b)K_1))$. This completes the proof of Theorem 1.5. \hfill $\Box$

\section*{Data availability statement}

My manuscript has no associated data.

\section*{Declaration of competing interest}

The author declares that he has no conflicts of interest to this work.

\medskip

\section*{Acknowledgments}

The author is very grateful to the anonymous reviewers for their valuable comments and corrections which result in an improvement of the original manuscript.
This work was supported by the Natural Science Foundation of Jiangsu Province (Grant No. BK20241949). Project ZR2023MA078 supported by Shandong Provincial
Natural Science Foundation.

\end{document}